\documentclass[11pt]{article}
\usepackage{latexsym,amsmath,amssymb}
\usepackage{color}
\usepackage[colorlinks]{hyperref}
\newtheorem{defi}{Definition}[section]
\newtheorem{lemma}[defi]{Lemma}

\newtheorem{theorem}[defi]{Theorem}

\newcommand{\subject}[1]{\begin{flushleft}\textbf{Mathematics  Subject Classification}: #1\end{flushleft}}
\newcommand{\keyword}[1]{\par\noindent \textbf{Keywords:} #1 }
\newcommand{\university}[1]{\\[3mm]{\small #1}}
\newcommand{\email}[1]{\begin{flushleft}\texttt{E-mail:\textcolor[rgb]{0.50,0.00,1.00}{#1}}\end{flushleft}}

\begin{document}

\title{Universally Koszul algebras defined by graphs}
\author{Rashid Zaare-Nahandi\university{Institute for Advanced Studies in Basic Sciences, Zanjan 45195, Iran}
 }
\date{}
\maketitle
\begin{abstract}
In this note, it is proved that a graphs is $(2K_2,P_4)$-free if and only if its edge ring is universally Koszul. Using properties of this family of graphs, we show that Universally Koszul algebras defined by graphs have linear minimal free resolution.
\keyword{Universally Koszul algebra, linear minimal free resolution, $(2K_2,P_4)$-free graph.}
\subject{13D02}
\end{abstract}
Let $K$ be a field and $R$ be a standard graded $K$-algebra. More precisely,  $R=\oplus_{i\in {\Bbb N}} R_i$
is of the form $K[X_1, \ldots, X_n]/I$ where $I$ is a homogeneous ideal. Let $S=K[X_1,\ldots,X_n]$ and $x_i=\overline{X}_i$ be image of $X_i$ in $R$.
The minimal $R$-free resolution of a graded $R$-module $M$ is said to be linear
if the matrices that represent the maps of the resolution have entries of degree 1.

\begin{defi}
The algebra $R$ is called Koszul if $K\simeq R/\langle x_1, \ldots, x_n\rangle$ has a linear minimal $R$-free
resolution. We say $R$ is universally Koszul if for each ideal $J$ of $R$ generated by elements of degree 1, $R/J$ has a linear $R$-free resolution.
\end{defi}

The following lemma is known \cite{conca1}.

\begin{lemma}
Let $R\simeq S/I$ be universally Koszul algebra. Then, $I$ is an ideal generated by homogeneous polynomials of degree two.
\end{lemma}

For a survey on Koszul algebras we refer the reader to the paper of Fr\"oberg~\cite{froberg1} . For generalities on universally Koszul algebras refer to papers of Conca~\cite{conca1} and~\cite{conca2}. The goal of this note is to classify the universally Koszul algebras defined by square-free quadratics, that is, algebras of the form $S/I$ such that $I$ is generated by square-free monomials of degree two.



Let $G$ be a graph on the vertex set $V = \{v_l,... ,v_n\}$. The edge ideal $I(G)$, associated to $G$, is the ideal
of $S$ generated by the set of all square-free monomials $X_iX_j$ provided that $v_i$ is adjacent to $v_j$ in $G$. The ideal $I(G)$ is called edge ideal of $G$ and the ring $S/I(G)$ is called edge ring of $G$.

\begin{defi}
Let $R_1=K[X_1,\ldots,X_n]/I$ and $R_2=K[Y_1,\ldots,Y_m]/J$ be two homogeneous $K$-algebras. Define the fiber product of $R_1$ and $R_2$ to be $R_1\circ R_2 = K[X_1,\ldots,X_n,Y_1,\ldots,Y_m]/(I+J+H)$ where $H=\langle X_iY_j, 1\leq i\leq n, 1\leq j\leq m\rangle$.
\end{defi}

Note that if $R_1$ is edge ring of a graph $G_1$ and $R_2$ is edge ring of $G_2$, then $R_1\circ R_2$ is edge ring of join graph of $G_1$ and $G_2$, which is a graph consisting of $G_1$ and $G_2$ and all vertices of $G_1$ are adjacent to all vertices of $G_2$.

\begin{lemma} {\rm \cite{conca2}}\label{conca}
\begin{itemize}
\item[(a)] A ring is universally Koszul if and only if $R[X]$ is universally Koszul. Where, $X$ is an indeterminate over $R$.
\item[(b)] The fiber product $R_1\circ R_2$ of algebras $R_1$ and $R_2$ is universally Koszul if and only if $R_1$
and $R_2$ are both universally Koszul.
\item[(c)] If $R$ is universally Koszul and $I$ is an ideal of $R$ generated by elements of degree 1,
then $R/I$ is universally Koszul.
\item[(d)] Let for any nonzero element $z$ of degree 1 in $R$, the algebra $R/(z)$ is universally Koszul and the ideal $(0:z)$ is generated by elements of degree 1. Then $R$ is universally Koszul.
\end{itemize}
\end{lemma}

Let $2K_2$ be disjoint union of two $K_2$ (complete graph with two vertices and one edge), and let $P_4$ be a path of length 3 (a graph with 4 vertices and 3 edges connecting them successively).  We say that a graph $G$ is $(2K_2,P_4)$-free if there is no any induced subgraph of $G$ isomorphic to $2K_2$ or $P_4$. A graph is called chordal if there is no any induced subgraph isomorphic to $C_l$ (cycle of length $l$) for $l>3$.

Now we state the main theorem of this note.

\begin{theorem} Let $R$ be the edge ideal of a graph $G$. The following conditions are equivalent:
\begin{itemize}
\item[(1)] $R$ is universally Koszul,
\item[(2)] The graph $G$ is $(2K_2,P_4)$-free.
\end{itemize}
\end{theorem}

{\it Proof.} First note that if we add some vertices with degree 0 to the graph $G$, then by Lemma~\ref{conca}(a), the ring $R(G)$ is universally Koszul if and only if the edge ring of the new graph is universally Koszul, therefore, we may assume that the graph $G$ does not have any vertex with degree 0. 

(1) $\to$ (2) Assume that $G$ is not $(2K_2,P_4)$-free. Then $G$ has an induced subgraph isomorphic to $P_4$ or $2K_2$. In this case, there are indeterminates $x_{i_1}, x_{i_2}, x_{i_3}, x_{i_4}$ such that at least one of the following conditions hold.
\begin{itemize}
\item[a)] $x_{i_1}x_{i_2}, x_{i_3}x_{i_4}\in I(G)$ and $x_{i_1}x_{i_3}, x_{i_1}x_{i_4}, x_{i_2}x_{i_3}, x_{i_2}x_{i_4}\not\in I(G)$.
\item[b)] $x_{i_1}x_{i_2}, x_{i_2}x_{i_3}, x_{i_3}x_{i_4}\in I(G)$ and  $x_{i_1}x_{i_3}, x_{i_1}x_{i_4}, x_{i_2}x_{i_4}\not\in I(G)$.
\end{itemize}
In both cases $x_{i_1}x_{i_4}(x_{i_2}+x_{i_3})=0$ in $R$ but $x_{i_1}(x_{i_2}+x_{i_3})\neq 0$ and $x_{i_4}(x_{i_2}+x_{i_3})\neq 0$. Therefore, $x_{i_1}x_{i_4}$ is in the minimal generating set of the ideal $0:(x_{i_2}+x_{i_3})$ and is not of degree 1. By Lemma~\ref{conca}(c,d), $R$ is not universally Koszul.

(2) $\to$ (1)  Assume that $G$ is $(2K_2,P_4)$-free. Let $N(v_i)$ be the set of all vertices in $G$ adjacent to $v_i$. The prof is by induction on number of vertices of $G$. If $n=|V(G)|=1$, then the assertion is clearly true. Assume that $n>1$. Let $A=\{v_1,\ldots,v_r\}$ be a maximal independent set of vertices in $G$ and $1\leq i\leq r$ be given. Then each vertex in $N(v_i)$ is adjacent to each vertex in $V(G)\setminus N(v_i)$. To see this, assume that $v\in N(v_i)$ and $w\in V(G)\setminus N(v_i)$.  Then $w\in A$ and there is $u\in V(G)$ such that $u\sim w$, or $w\not\in A$ and there is $u\in A$ such that $u\sim w$. In both cases, $v_i$ is adjacent to none of $u$ and $w$. Then $v$ must be adjacent to both $u$ and $w$, because in other case we have an induced subgraph of $G$ isomorphic to $2K_2$ or $P_4$. Therefore, $G$ is join of two induced graphs on $N(v_i)$ and $V(G)\setminus N(v_i)$ and by induction hypothesis, both graphs are universally Koszul and by Lemma~\ref{conca}(b), $G$ is universally Koszul. \hfill $\Box$

Let $G$ be a graph. Then $G$ is $(2K_2,P_4)$-free graph if and only if $\overline{G}$, the complement of $G$, is $(C_4,P_4)$-free, where $C_4$ is the cycle of length 4. Therefore, in $\overline{G}$ there is no any induced cycle of length 4 or greater and it is a chordal graph.

We use the following theorem of Fr\"oberg in reformulated form to state the final step of this note.
\begin{theorem} {\rm [4]} Let $I = I(G)$ be edge ideal of a graph $G$. Then, $S/I$ has linear minimal $S$-free resolution if and only if $\overline{G}$ is a chordal graph.
\end{theorem}

\begin{theorem}
Let $R=S/I$ be a universally Koszul algebra where $I=I(G)$ is edge ideal of a graph $G$. Then $R$ has linear minimal free resolution as $S$-module.
\end{theorem}

{\it Proof.} The proof is clear using the previous Theorem and the fact that $\overline{G}$ is chordal. \hfill $\Box$

\email{rashidzn@iasbs.ac.ir}
\end{document}